\documentclass[MSNbibl,number,citesort,seceqn,dvips]{arxbj}
\usepackage{upgreek}

\aid{0}
\volume{19}
\issue{5A}
\pubyear{2013}
\firstpage{1818}
\lastpage{1838}
\doi{10.3150/12-BEJ431} 

\makeatletter
\newtheorem{CORR}{Corollary}[section]
\newtheorem{CONJ}{Conjecture}[section]

\newtheorem{lemm}{Lemma}[section]
\newtheorem{prop}{Proposition}[section]
\newtheorem{THEO}{Theorem}[section]
\theoremstyle{remark}

\newremark{REMS}{Remarks}[section]

\newcommand{\ba}{b_\alpha}
\newcommand{\bz}{{\bar z}}
\newcommand{\ca}{c_\alpha}
\newcommand{\CC}{{\mathbb{C}}}
\newcommand{\cP}{{\mathcal{P}}}
\newcommand{\Ea}{E_\alpha}
\newcommand{\EE}{{\mathbb{E}}}
\newcommand{\fa}{f_\alpha}
\newcommand{\ga}{g_\alpha}
\newcommand{\ha}{h_\alpha}
\newcommand{\hba}{{{\tilde b}_\alpha}}
\newcommand{\hta}{{{\tilde a}_\alpha}}
\newcommand{\hca}{{{\tilde c}_\alpha}}
\newcommand{\lbd}{\lambda}
\newcommand{\Ma}{M_\alpha}
\newcommand{\rl}{{\mathbb{R}}}
\newcommand{\Ta}{T_\alpha}
\newcommand{\Un}{\mathbf{1}}
\newcommand{\Va}{V_\alpha}
\newcommand{\Xa}{X_\alpha}
\newcommand{\Ya}{Y_\alpha}
\newcommand{\Wa}{W_\alpha}
\newcommand{\Za}{Z_\alpha}

\newcommand{\elaw}{\,\stackrel{d}{=}\,}
\makeatother

\begin{document}
\begin{frontmatter}

\title{Further examples of GGC and HCM densities}

\runtitle{GGC and HCM densities}

\begin{aug}
\author[1]{\fnms{Wissem} \snm{Jedidi}\thanksref{1}\ead[label=e1]{wissem.jedidi@ksu.edu.sa}}
\and
\author[2]{\fnms{Thomas} \snm{Simon}\corref{}\thanksref{2}\ead[label=e2]{simon@math.univ-lille1.fr}}
\runauthor{W. Jedidi and T. Simon} 
\address[1]{Department of Statistics and Operations Research, College of Science, King
Saud University, P.O. Box 2455, Riyadh 11451, Saudi Arabia.
\printead{e1}}
\address[2]{Laboratoire Paul Painlev\'e, Universit\'e Lille 1, Cit\'e
Scientifique, 59655 Villeneuve d'Ascq Cedex, France. \printead{e2}}
\end{aug}

\received{\smonth{4} \syear{2011}}
\revised{\smonth{1} \syear{2012}}

%
\begin{abstract}
We display several examples of generalized gamma convoluted and
hyperbolically completely monotone random variables related to positive
$\alpha$-stable laws. We also obtain new factorizations for the latter,
refining Kanter's and Pestana--Shanbhag--Sreehari's. These results give
stronger credit to Bondesson's hypothesis that positive $\alpha$-stable
densities are hyperbolically completely monotone whenever
$\alpha\le1/2.$
\end{abstract}

%
\begin{keyword}
\kwd{generalized Gamma convolution}
\kwd{hyperbolically completely monotone}
\kwd{hyperbolically monotone}
\kwd{positive stable density}\vspace*{12pt}
\end{keyword}

\end{frontmatter}

\section{Introduction}\label{s1}

A positive random variable $X$ is called a Generalized Gamma
Convolution (GGC) if its Laplace transform reads
\begin{equation}
\label{grege} \EE[\mathrm{e}^{-\lbd X}] = \exp- \biggl[a\lbd+
\int_0^\infty(1- \mathrm{e}^{-\lbd x}) \frac{\varphi(x)}{x}
\,\mathrm{d}x\biggr],\qquad \lbd\ge0,
\end{equation}
where $a\ge0$ and $\varphi$ is a completely monotone (CM) function over
$(0, +\infty).$ The denomination comes from the fact that the above
class can be identified as the closure for weak convergence of finite
convolutions of Gamma distributions. We refer to \cite{B} and
\cite{SVH} for comprehensive monographs on such random variables. From
their definition, GGC random variables are self-decomposable (SD) hence
infinitely divisible (ID), absolutely continuous and unimodal -- see,
for example, \cite{S} for the proofs of the latter properties. We also
see from (\ref{grege}) that GGC random variables are characterized up
to translation by the positive Radon measure on $(0,+\infty)$ uniquely
associated to the CM function $\varphi$ by Bernstein's theorem, which
is called the Thorin measure of $X$ and whose total mass,
$\varphi(0+),$ might be infinite. As an illustration of this
characterization, Theorem 4.1.4\vadjust{\eject} in \cite{B} shows that the density of
$X$ vanishes in $a+$ if $\varphi(0+) > 1,$ whereas it is infinite in
$a+$ if $\varphi(0+) < 1.$ We refer to \cite{JRY} for a recent survey
on GGC variables having a finite Thorin measure, dealing in particular
with their Wiener--Gamma representations and their relations with
Dirichlet processes.

A positive random variable $X$ is said to be hyperbolically completely
monotone (HCM) if it has a density $f$ on $(0, +\infty)$ such that for
every $u > 0$ the function
\[
H_u(w) = f(uv)f(u/v),\qquad w = v + 1/v \ge2,
\]
is CM in the variable $w$ (it is easy to see that $H_u$ is always a
function of $w$). In general, a function $f \dvtx
(0,+\infty)\to(0,+\infty)$ is said to be HCM when the above CM property
holds for $H_u,$ and this extended definition will be important in the
sequel. HCM densities turn out to be characterized as pointwise limits
of densities of the form
\begin{equation}
\label{dens} x \mapsto C x^{\beta- 1}\prod_{i = 1}^N (x +
y_i)^{-\gamma_i},
\end{equation}
where all above parameters are positive -- see Sections 5.2 and 5.3 in
\cite{B}. This characterization yields many explicit examples of GGC
random variables, since it is also true that HCM random variables are
GGC -- see Theorem 5.2.1 in \cite{B}. Actually, HCM variables appear as
a kind of center for GGC in view of Theorem 6.2.1 in \cite{B} which
states that the independent product or quotient of a GGC by a HCM
variable is still a GGC. The HCM class is also stable by independent
multiplication and power transformations of absolute value greater than
one. We refer to \cite{B} for many other properties of HCM densities
and functions.

The HCM property is connected to log-concavity in the following way. A
positive random variable $X$ is said to be hyperbolically monotone (HM)
if it has a density $f$ on $(0, +\infty)$ such that the above function
$H_u$ is nonincreasing in the variable $w$. Similarly as above, one can
extend the HM property to all positive functions on $(0, +\infty).$
Obviously, HCM is a subclass of HM. It is easy to see -- see \cite{B},
pages 101--102 -- that $X$ is HM iff its density $f$ is such that
$t\mapsto f(\mathrm{e}^t)$ is log-concave on $\rl$. This shows that $f$
is a.e. differentiable with $x \mapsto x f'(x)/f(x)$ a nonincreasing
function, so that $f'$ has at most one change of sign and $X$ is
unimodal. The main theorem in \cite{CT} shows that HM variables are
actually multiplicatively strong unimodal, viz. their independent
product with any unimodal random variable is unimodal. From the
log-concavity characterization, the HM property is stable by power
transformation of any value, and this entails that the inclusion HCM
$\subset$ HM is strict: if $L\sim \operatorname{Exp} (1),$ then
$\sqrt{L}$ is HM but not ID, hence not HCM. From Pr\'ekopa's theorem,
the HM property is also stable by independent
multiplication.\looseness=1

For a positive random variable with density, the standard way to derive
the GGC property is to read it from the Laplace transform. For example,
it is straightforward to see that positive $\alpha$-stable variables
are all GGC -- see Example 3.2.1 in \cite{B}. On the other hand, it is
easier to study the HCM property from features on the\vadjust{\eject} density itself
and Laplace transforms are barely helpful. As an illustration of this,
we show without much effort in Section~\ref{s4} of the present paper
that the quotient of two positive $\alpha$-stable variables, whose
density is explicit, is HCM iff $\alpha\le1/2.$ Problems become usually
more intricate when one searches for GGC without explicit Laplace
transform or for HCM without closed expression for the
density.\looseness=1

In 1981, Bondesson raised the conjecture that positive $\alpha$-stable
variables should be HCM iff $\alpha\le1/2$ and this very hard problem
(quoting his own recent words, see Remark~3 in~\cite{B09}) is still
unsolved except in the easy case when $\alpha$ is the reciprocal of an
integer -- see Example 5.6.2 in \cite{B}. Notice, in passing, that the
validity of this conjecture is erroneously taken for granted in~\cite{JRY},
page 361. We refer to \cite{B81}, pages 54--55, \cite{B}, pages 88--89 and also to the manuscript~\cite{Bo99},
for several reasons, partly numerical,
supporting this hypothesis. Let us also mention the main theorem of
\cite{TS}, which states that positive $\alpha$-stable random variables
are HM iff $\alpha\le 1/2.$ Actually, it follows easily from the proofs
of Lemmas 1 and 2 in \cite{TS} that the $p$th power of a positive
$(p/n)$-stable variable is HCM for any integers $p, n\ge2$ such that
$p/n \le1/2.$

In the present paper, we will present several examples of GGC and HCM
densities related to the above conjecture. In Section~\ref{s2}, we
combine the main results of \cite{RVY} and \cite{TS} to show the GGC
property for a large family of negative powers of $\alpha$-stable
variables with $\alpha\le 1/4$. This family is actually a bit larger
than the one which would be obtained from the validity of Bondesson's
hypothesis. The more difficult case $\alpha\in(1/4, 1/2]$ is also
studied, with a partial result. In Section~\ref{s3}, we use Kanter's
and Pestana--Shanbhag--Sreehari's factorizations to show that a large
class of positive powers of $\alpha$-stable variables is the product of
an HCM variable and an ID variable. The latter turns out to be always a
mixture of exponentials (ME), hence very close to a GGC. Along the way,
we also obtain an independent proof of Pestana--Shanbhag--Sreehari's
factorization. In Section~\ref{s4}, we show the aforementioned HCM
result for the quotient of two stable variables, and a similar
characterization for Mittag--Leffler variables. Not surprisingly, both
yield the same boundary parameter $\alpha=1/2.$

The results presented in Sections \ref{s2} and \ref{s3} are probably
not optimal and at the end of Section~\ref{s3} we state another
conjecture, where the power exponent $\alpha/(1-\alpha)$ appears
naturally. We also hope that the different tools and methods presented
here will be helpful to tackle Bondesson's conjecture more deeply, even
though we have tried to exploit them to their full extent.

\subsection*{Notations}

We will consider real random variables $X$ having a density always
denoted by $f_X$, unless explicitly stated. For the sake of brevity, we
will use slightly incorrect expressions like ``GGC variable'' or ``HCM
variable'' and sometimes even delete the word ``variable'' (as was
actually already done in the present introduction). We will also set
``positive (negative) $\alpha$-stable power'' for ``positive (negative)
power transformation of a positive $\alpha$-stable random variable''.

\section{\texorpdfstring{Negative $\alpha$-stable powers and the GGC property}
{Negative alpha-stable powers and the GGC property}}\label{s2}

\subsection{Some consequences of the HM property}\label{s2.1}

Let $\Za$ be a positive $\alpha$-stable random variable --
$\alpha\in(0,1)$ -- with density function $\fa$ normalized such that
\[
\int_0^\infty \mathrm{e}^{-\lbd t} \fa(t) \,\mathrm{d}t =
\EE[\mathrm{e}^{-\lbd\Za}]= \mathrm{e}^{-\lbd^\alpha},\qquad \lbd\ge0.
\]
In the remainder of this paper, we will use the notation $\beta= 1
-\alpha.$ We will also set $Z_1 = 1$ by continuity. Recall that when
$\alpha= 1/2,$ our normalization yields
\begin{equation}
\label{half} f_{1/2}(x) = \frac{1}{2\sqrt{\uppi}
x^{3/2}}\mathrm{e}^{-1/4x}\Un _{\{ x > 0\}}\cdot
\end{equation}
Kanter's factorization -- see Corollary 4.1 in \cite{K} -- reads
\begin{equation}
\label{Kant} \Za\elaw L^{-\beta/\alpha}\times b_\alpha^{-1/\alpha}(U),
\end{equation}
where $L\sim\operatorname{Exp} (1),$ $U\sim
\operatorname{Unif}(0,\uppi)$ independent of $L$, and
\[
b_\alpha(u) = \bigl(\sin u/\sin(\alpha u)\bigr)^\alpha\bigl(\sin
u/\sin(\beta u)\bigr)^\beta,\qquad u\in (0,\uppi),
\]
is a bounded, decreasing and concave function -- see Lemma 1 in
\cite{TSS}. Observe that when $\alpha= 1/2,$ Kanter's factorization is
a particular instance of the so-called Beta--Gamma algebra -- see, for
example, \cite{B}, pages 13--14. Indeed, one has
\[
4 b_{1/2}^{-2} (U) = \cos^{-2}(U/2)\elaw\operatorname{Beta}^{-1} (1/2,
1/2)
\]
and (\ref{half}) entails $4 Z_{1/2} \elaw\operatorname{Gamma}^{-1}
(1/2, 1),$ so that (\ref{Kant}) amounts when $\alpha= 1/2$ to
\[
\operatorname{Gamma} (1/2, 1) \elaw\operatorname{Beta} (1/2, 1/2)
\times\operatorname{Gamma} (1,1).
\]
Put together with Shanbhag--Sreehari's classical factorization of the
exponential law -- see, for example, Exercise 29.16 in \cite{S},
Kanter's factorization also shows that for every
$\gamma\ge\alpha/\beta$ the random variable $\Za^{-\gamma}$ is ME, viz.
there exists a positive random variable $U_{\alpha, \gamma}$ such that
\begin{equation}
\label{FactPSS} \Za^{-\gamma} \elaw L \times U_{\alpha, \gamma}.
\end{equation}
See, for example, Section~51.1 in \cite{S} for more material on ME
random variables. With the help of the HM property, one has the
following reinforcement.

\begin{prop}
\label{basic} For every $\gamma> 0,$ the random variable
$\Za^{-\gamma}$ is ID (with a CM density) iff $\gamma\ge \alpha/\beta.$
Moreover, $\Za^{-\gamma}$ is SD for every $\alpha\le 1/2$ and every
$\gamma\ge\alpha/\beta.$
\end{prop}

\begin{pf} The factorization (\ref{FactPSS}) together with Theorem 51.6 and
Proposition 51.8 in \cite{S} show that $\Za^{-\gamma}$ is ID with a CM
density if $\gamma\ge\alpha/\beta.$ On the other hand, a change of
variable and Linnik's asymptotic expansion -- see, for example, (14.35)
in \cite{S} -- yield
\[
x^{\alpha/\beta\gamma} \log f_{\Za^{-\gamma}}(x) \to\kappa_{\alpha,
\gamma } \in (-\infty, 0)
\]
as $x\to\infty$ for every $\gamma> 0.$ Hence, if $\gamma<
\alpha/\beta,$ Theorem 26.1 in \cite{S} -- see also Exercise 29.10
therein -- entails that $\Za^{-\gamma}$ is not ID. When $\alpha\le1/2,$
the main result in \cite{TS} shows that $\Za$ is HM, so that
$\log(\Za^{-\gamma})$ has a log-concave density. If in addition
$\gamma\ge\alpha/\beta,$ we have just observed that $\Za^{-\gamma}$ has
a CM density which is hence decreasing and log-convex. This entails
that when $\alpha\le1/2$ and $\gamma\ge\alpha/\beta,$ the random
variable $\Za^{-\gamma}$ belongs to the class mentioned in \cite{B},
Remark~VI, page~28, and is SD.
\end{pf}

The main result of this section shows that the SD property for
$\Za^{-\gamma}$ can be refined into GGC, in some cases. The proof also
relies on the HM property.

\begin{THEO}
\label{4a} The random variable $\Za^{-\gamma}$ is GGC for any
$\alpha\in(0,1/4], \gamma\ge4\alpha.$
\end{THEO}

\begin{pf}
Fix $\alpha\in(0,1/4], \gamma\ge4\alpha$ and set $\delta=
2\alpha/\gamma\in(0, 1/2].$ Bochner's subordination for stable
subordinators -- see, for example, Example 30.5 in \cite{S} -- yields
the identity
\[
\Za^{-2\alpha} \elaw\ca(Z_{1/2}^{-1}\times Z_{2\alpha}^{-2\alpha})
\]
for some purposeless constant $\ca> 0.$ We hence need to show the GGC
property for the random variable
\[
\bigl((4Z_{1/2})^{-1} \times Z_{2\alpha}^{2\alpha}\bigr)^{1/\delta},
\]
whose density is in view of (\ref{half}), the multiplicative
convolution formula, and a series of standard changes of variable,
expressed as
\[
\frac{\delta x^{\delta-3/2}}{\alpha\sqrt{2\uppi}}\int_0^\infty
\mathrm{e}^{-x^\delta y} f_{2\alpha}(y^{1/2\alpha}) y^{1/2\alpha- 1/2}
\,\mathrm{d}y.
\]
Observe that since
\[
\int_0^\infty f_{2\alpha}(y^{1/2\alpha}) y^{1/2\alpha- 1/2}
\,\mathrm{d}y = 2\alpha\int _0^\infty f_{2\alpha}(z) z^{\alpha}
\,\mathrm{d}z = \frac{2\alpha\sqrt{\uppi}}{\Gamma(1-\alpha)} < +\infty
\]
(see, e.g., (25.5) in \cite{S} for the second equality), the function
$K_\alpha f_{2\alpha}(y^{1/2\alpha}) y^{1/2\alpha- 1/2}$ is a
probability density on $\rl^+,$ where we have set $K_\alpha=
\Gamma(1-\alpha)/2\alpha\sqrt{\uppi}.$ Denoting by $\Xa$ the
corresponding random variable, we have to show that
\[
x \mapsto K_{\alpha, \delta}x^{\delta-
3/2}\EE[\mathrm{e}^{-x^\delta\Xa}]
\]
is the density of a GGC, with $K_{\alpha, \delta} =
\sqrt{2}\delta/\Gamma(1-\alpha).$ Since $\delta\le1/2 < 3/2,$ we see
from Theorem~6.2 in \cite{BB} (and the Remark 6.1 thereafter) that this
will be done as soon as
\[
\EE[\mathrm{e}^{-x^\delta\Xa}] = \EE\bigl[\mathrm{e}^{- x
(Z_\delta\times\Xa^{1/\delta})}\bigr]
\]
is, up to normalization, the density of a GGC. We will now obtain this
property with the help of Theorem 2 in \cite{RVY}. On the one hand, it
is easy to see that all negative moments of $Z_\delta$ and $\Xa$ are
finite, so that the density of $Z_\delta\times\Xa^{1/\delta}$ fulfils
(1.1) in \cite{RVY}. On the other hand, the main result of \cite{TS}
entails that $Z_{2\alpha}$ is HM because $2\alpha\le1/2,$ so that the
function
\[
t \mapsto f_{2\alpha}(\mathrm{e}^{t/2\alpha}) \mathrm{e}^{t/2\alpha-
t/2}
\]
is log-concave and $\Xa$ is HM as well. Also, $Z_\delta$ is HM because
$\delta\le1/2.$ Since the HM property is stable by independent
multiplication, this shows that $Z_\delta\times\Xa^{1/\delta}$ is HM,
in other words that it belongs to the class ${\mathcal{C}}$ defined in
\cite{RVY}, page 183, and we can apply Theorem 2 therein to conclude
the proof.
\end{pf}

\begin{REMS}
\label{4ar}
(a) From (14.30) in \cite{S} and a change of variable, one
has
\[
\sup\Bigl\{ u; \lim_{x \to0} f_{\Za^{-\gamma}}(x)/x^{u-1} = 0\Bigr\} =
\alpha /\gamma
\]
for every $\gamma> 0,$ so that (3.1.4) in \cite{B} shows that under the
assumptions of Theorem \ref{4a}, the GGC random variable
$\Za^{-\gamma}$ has a finite Thorin measure whose total mass is
$\alpha/\gamma.$ In other words, there exists a nonnegative random
variable $G_{\alpha,\gamma}$ such that
\[
\EE[\mathrm{e}^{-\lbd\Za^{-\gamma}}] = \exp-\biggl[(\alpha/\gamma)\int
_0^\infty(1- \mathrm{e}^{-\lbd x})\EE[\mathrm{e}^{-x
G_{\alpha,\gamma}}]\frac{\mathrm{d}x}{x} \biggr],\qquad \lbd\ge0.
\]
It would be interesting to get more properties of the random variables
$G_{\alpha,\gamma}.$

(b) It is easily seen that the above proof remains unchanged (and is
even shorter) if we take $\delta= 1$ viz. $\gamma= 2\alpha,$ so that
$\Za^{-2\alpha}$ is GGC as well for any $\alpha\in(0,1/4].$ In view of
Theorem \ref{4a} and the general conjecture made in \cite{B09},
Remark~3(ii), it is plausible that $\Za^{-\gamma}$ is GGC for any
$\alpha\in (0,1/4], \gamma\ge2\alpha.$
\end{REMS}

\subsection{A certain family of densities on $\rl^+$ and a partial result}\label{s2.2}

A drawback of Theorem \ref{4a} is that it only covers the range
$\alpha\in(0, 1/4].$ Indeed, with the same subordination method one
should expect to handle the range $\alpha\in(1/4, 1/2]$ as well.
Motivated by the key-properties (2.20) and (2.23) in the proof of
Theorem 2 in \cite{RVY}, let us define the class $\cP$ of probability
densities $f$ on $(0, +\infty)$ satisfying
\begin{equation}
\label{CPP} f(x)f(c/x) \ge f(1/x)f(cx)
\end{equation}
for all $x, c > 0$ such that $(x-1)(c-1) \ge0.$ With an abuse of
notation, we shall say that a random variable $X$ with density $f$
belongs to $\cP$ if $f\in\cP.$ If $X\in\cP,$ then it is easy to see
that $X^\gamma\in\cP$ for any $\gamma\neq0.$ Besides, it follows from
\cite{RVY}, pages 187--188, that HM $\subset\cP.$ Notice also that $\cP
\not \subset$ HM, as the following example shows. Consider the
independent quotient $\Ta= (\Za/\Za)^\alpha$ which has an explicit
density $\ga$ given by
\[
\ga(x) = \frac{\sin\uppi\alpha}{\uppi\alpha(x^2 + 2 \cos(\uppi\alpha) x
+1)}
\]
(see, e.g., Exercise 4.21(3) in \cite{CY}). A computation yields
\[
\frac{(\uppi\alpha)^2}{\sin^2
\uppi\alpha}\biggl(\frac{1}{\ga(x)\ga(c/x)} - \frac
{1}{\ga(cx)\ga(1/x)}\biggr)= (1- c^2) (x - 1/x)(x + 1/x + 2
\cos\uppi\alpha),
\]
which is clearly nonpositive whenever $(x-1)(c-1) \ge0,$ so that
$\Ta\in\cP$ for any $\alpha\in(0,1).$ However, it is easy to show --
see the proof of Corollary \ref{co4.1} below -- that $\Ta$ is HM iff
$\alpha\le1/2.$ However, we know from (ix), page 68 in \cite{B} that
$\Ta$ is never HCM since $\ga$ has two poles
$\mathrm{e}^{\mathrm{i}\uppi\alpha}$ and
$\mathrm{e}^{-\mathrm{i}\uppi\alpha }$ in $\CC\setminus(-\infty, 0].$
Notice also that the variable $T_{1/2}$ is SD but not GGC -- see
\cite{Di} and the references therein. We will come back to this example
in Section~\ref{s4}. The following proposition makes the relationship
between HM and $\cP$ more precise.

\begin{prop}
\label{LM} For any nonnegative random variable $X$ having a density,
one has
\[
X\mbox{ is }\mathrm{HM}\quad\Longleftrightarrow\quad
cX\in\cP\qquad\forall c
> 0.
\]
\end{prop}

\begin{pf} The direct part is easy since $cX$ is HM for any $c > 0$
whenever $X$ is HM. For the indirect part, setting $g_X(t) = \log
f_X(\mathrm{e}^t)$ for any $t \in\rl,$ the fact that $cX\in\cP$ for any
$c > 0$ shows that for any $-\infty< a\le b\le c\le d < +\infty$ with
$b+c = a+d,$ one has $g_X(b) + g_X(c) \ge g_X(a) + g_X(d),$ so that
$g_X$ is concave.
\end{pf}

Together with the above example, this proposition entails that $\cP$ is
not stable by multiplication with positive constants. Since on the
other hand $\cP$ is clearly stable under weak convergence and since any
positive constant can be approximated by a sequence of truncated
gaussian variables which all belong to HM $\subset\cP$, the instability
of $\cP$ w.r.t. constant multiplication entails that $\cP$ is not --
contrary to HM -- stable by independent multiplication either, viz.
there exist independent $X, Y \in\cP$ such that $X\times Y \notin \cP.$

Let now $Y$ be a nonnegative random variable with a density of the form
$\kappa x^{-a}\EE[\mathrm{e}^{-xX}]$ for some $a\ge0$ and a nonnegative
random variable $X$ with finite negative moments such that
$c^{-1}X\in\cP$ for some $c>0.$ A perusal of the proof of Theorem 2 in
\cite{RVY} -- see especially (2.10), (2.20) and (2.23) therein -- shows,
together with Theorem 6.2 in \cite{BB}, that $\varphi_Y(x) =
\EE[\mathrm{e}^{-x Y}]$ is such that $H_c(w) = \varphi_Y(c
v)\varphi_Y(c/v)$ is CM in the variable $w = v + 1/v.$ On the other
hand, it is possible to show the following intrinsic property of
positive stable densities.

\begin{THEO}
\label{cphcm} For every $\alpha\in(0,1)$, there exists $c_\alpha\ge0$
such that $c\Za\in\cP\Leftrightarrow c\ge c_\alpha.$
\end{THEO}

Though it has independent interest, we prefer not giving the proof of
this theorem since it is quite long, relying on the single intersection
property for $\Za$ -- see Theorem 4.1 in \cite{K}, an extended Yamazato
property for $\fa$ which is displayed in (1.4) in \cite{TS} and the
discussion thereafter, and a detailed analysis. Notice from Proposition
\ref{LM} and the main result in \cite{TS} that $c_\alpha= 0$ for any
$\alpha \le1/2,$ and that necessarily $c_\alpha> 0$ when $\alpha> 1/2.$
Theorem \ref{cphcm} and a painless adaptation of the proof of Theorem
\ref{4a} entail the following property of the variable
$\Za^{-2\alpha}.$

\begin{CORR}
\label{party} For every $\alpha\in(1/4, 1/2],$ there exists $\hca> 0$
such that for every $c \in[0, \hca],$ the function $H^\alpha_c(w) =
\varphi_\alpha(c v)\varphi_\alpha(c/v)$ is CM in the variable $w = v +
1/v,$ where $\varphi_\alpha(x) = \EE[\mathrm{e}^{-x \Za^{-2\alpha}}],
x\ge0.$
\end{CORR}

If we could show that $\hca= +\infty,$ then Theorem 6.1.1 in \cite{B}
would entail that $\Za^{-2\alpha}$ is GGC for every $\alpha\in(1/4,
1/2].$ Notice from Proposition~\ref{basic} that when $\alpha> 1/2$ the
random variable $\Za^{-2\alpha}$ is not ID (since then $2\alpha<
\alpha/\beta$), hence not GGC. From Corollary~\ref{party} and the above
Remark~\ref{4ar}(b), it is very plausible that $\Za^{-\gamma}$ is GGC
for any $\alpha\le1/2, \gamma\ge 2\alpha.$ See also the general
Conjecture~\ref{c3} raised at the end of the next section.

\section{On the infinite divisibility of Kanter's random variable}\label{s3}

In this section, we will deal with positive $\alpha$-stable powers.
From the point of view of the factorization (\ref{Kant}), we need to
study negative powers of $L$ and positive powers of the random variable
$b_\alpha^{-1/\alpha}(U),$ which will be referred to as Kanter's
variable subsequently. The latter plays an important role in simulation
-- see \cite{De} where it is called Zolotarev's variable, although the
original computation leading to (\ref{Kant}) is due to Chernine and
Ibragimov as explained in \cite{K}. It is interesting to remark that
Kanter's variable also appears explicitly in the context of free stable
laws -- see \cite{D}, page 138 and the references therein.

Negative powers of $L$ are not completely well understood from the
point of view of infinite divisibility. It follows from (iv) in
\cite{B} that $L^s$ is HCM for every $s\le-1,$ but it is not even known
whether $L^s$ is ID or not for $s\in(-1,0)$ -- see \cite{SVH}, page
521. Here, we will rather focus on positive powers of Kanter's
variable. Let us first notice that the above factorization
(\ref{FactPSS}) was also observed in \cite{PSS}, Theorem 1, where it is
actually shown that $U_{\alpha, \gamma} = \exp(-W_{\alpha, \gamma})$
for some ID random variable $W_{\alpha, \gamma}$. In this section, we
will investigate (\ref{FactPSS}) more thoroughly and give finer
properties of the random variable
\[
\Va= U_{\alpha, \alpha/\beta}^{-1} = b_\alpha^{-1/\beta}(U).
\]
We could actually consider any positive power of Kanter's random
variable, but the latter choice is more convenient for our purposes
because of the identities
\begin{equation}
\label{GGCKant} \Za^{\alpha s/\beta} \elaw L^{-s} \times\Va^s
\end{equation}
for every $s \in\rl.$ Our purpose is three-fold. First, we provide an
alternative proof of Theorem 1 in \cite{PSS}, with the improvement that
each positive power $\Va^s$ (in particular, Kanter's variable itself)
is ID with a log-convex density. Second, we show that all $\Va^s$ are
actually positively translated ME's. Third, we study in some detail the
case $\alpha= 1/2$ and propose a general conjecture which is, in some
sense, a reinforcement of Bondesson's.

\subsection{Another proof of Pestana--Shanbhag--Sreehari's
factorization}\label{s3.1}

Let us consider the random variable
\[
\Wa= \log(\Va) = -(1/\beta) \log(\ba(U))
\]
and observe that $W_{\alpha, \gamma} \elaw\gamma_\alpha^{-1}\Wa+ \log
(Z_{\gamma_\alpha})$ for every $\gamma\ge\alpha/\beta,$ with the
notation $\gamma_\alpha= \alpha/\beta\gamma.$ Since
$\log(Z_{\gamma_\alpha})$ is ID -- see, for example, Exercise 29.16 and
Proposition 15.5 in \cite{S}, Theorem 1 in \cite{PSS} follows as soon
as $\Wa$ is ID. In view of Theorem 51.2 in \cite{S} and the fact
(obvious from the definition of $\ba$) that the support of $\Wa$ is
unbounded on the right, this is a consequence of the following, which
we prove independently of \cite{PSS}.

\begin{THEO}
\label{LCX1} The density of $\Wa$ is log-convex.
\end{THEO}

This theorem entails that all positive powers of Kanter's random
variable $\ba^{-1/\alpha}(U)$ are ID, as shown in the next corollary.
In particular, $\Za^\gamma$ is the product of a HCM random variable and
an ID random variable for every $\gamma\ge\alpha/\beta.$

\begin{CORR}
\label{LCX2} For every $s > 0,$ the density of $\Va^s$ is decreasing
and log-convex.
\end{CORR}

\begin{pf}
Suppose first that $s = 1.$ Since the support of $\Wa$ is
unbounded on the right, the log-convexity of its density entails that
it is also decreasing, so that the function $\ga(t) = f_{\Va}
(\mathrm{e}^t)$ is also decreasing and log-convex. Since $\log x$ is
increasing and concave, this shows that $f_{\Va} (x) = \ga(\log x)$ is
decreasing and log-convex. The general case $s > 0$ follows analogously
in considering the variable $s\Wa= \log(\Va^s).$
\end{pf}

The proof of Theorem \ref{LCX1} relies on the following lemma.

\begin{lemm}
\label{trigo} The function $\ha= \ba''\ba/(\ba')^2$ is increasing on
$(0,\uppi)$.
\end{lemm}

\begin{pf}
First, observe that $b_{1/2} (u) = 2 \cos(u/2),$ so that
$h_{1/2} (u) = -\cot^2(u/2),$ an increasing function on $(0, \uppi)$.
In the general case, the proof is more involved. Since $b_\alpha=
b_\beta,$ it is enough to consider the case $\alpha< 1/2.$ Set
$A_\gamma(u) = \gamma \cot(\gamma u) - \cot(u)$ for every
$\gamma\in(0,1),$ $f = \alpha A_\alpha+ \beta A_\beta$ and $g = f' -
f^2.$ One has $\ba' = - f\ba$ and $\ba'' = - g\ba,$ so that $\ha= 1 +
(1/f)'$ and we need to show that $1/f$ is strictly convex, in other
words that
\begin{equation}
\label{cvx} 2(f')^2 - f f'' = 2 g f' - fg' > 0.
\end{equation}
It is shown in Lemma 1 of \cite{TSS} that $g = \alpha\beta+ h + k,$
with the further notations $h = \alpha\beta(A_\alpha- A_\beta)^2$ and
$k = 2(A_\alpha- A_\beta )(\beta A_\beta- \alpha A_\alpha).$ On the
other hand, the Eulerian formula
\[
\uppi\cot(\uppi z) = \frac{1}{z} + 2z \sum_{n\ge1} \frac{1}{z^2 - n^2}
\]
shows that
\[
A_\gamma(\uppi z) = \frac{2 (1 - \gamma^2)z}{\uppi}\sum_{n\ge1}\frac
{n^2}{(n^2 -z^2)(n^2 - \gamma^2z^2)}
\]
is a strictly absolutely monotonic function on $(0,1)$ (i.e., all its
derivatives are positive) for every $\gamma\in(0,1),$ so that $f =
\alpha A_\alpha+ \beta A_\beta$ is absolutely monotonic on $(0,
\uppi),$ too. In particular, since $\alpha\beta> 0,$ (\ref{cvx}) holds
if
\[
(2 h f' - fh') + (2k f' - fk') \ge0.
\]
A further computation entails $2 h f' - fh' = 2 (A_\alpha-A_\beta)
(A_\beta'A_\alpha- A_\beta A_\alpha')$ and $2 h f' - fh' = 2((\beta
A_\beta- \alpha A_\alpha)-2\alpha\beta(A_\alpha-A_\beta))
(A_\beta'A_\alpha- A_\beta A_\alpha'),$ so that we need to prove
\[
2(A_\beta'A_\alpha- A_\beta A_\alpha') \bigl((\alpha^2 + \beta^2)
(A_\alpha- A_\beta) + (\beta A_\beta - \alpha A_\alpha)\bigr) \ge0.
\]
The above Eulerian formula entails readily that $A_\alpha- A_\beta$ and
$\beta A_\beta- \alpha A_\alpha$ are positive (actually, absolutely
monotonic) functions on $(0,\uppi)$ and we are finally reduced to prove
\begin{equation}
\label{Final} A_\beta'A_\alpha- A_\alpha' A_\beta\ge0.
\end{equation}
We found no direct argument for the nonnegativity of the above
Wronskian. Writing
\[
\biggl(\frac{A_\beta}{A_\alpha}\biggr)(u) = \frac{\beta\cot(\beta u) -
\cot (u)}{\alpha \cot(\alpha u) - \cot(u)} = \biggl(\frac{\sin(\alpha
u)}{\sin(\beta u)}\biggr) \biggl( \frac{\beta\cos(\beta u)\sin(u) -
\cos(u)\sin(\beta u)}{\alpha\cos(\alpha u)\sin (u) - \cos(u)\sin(\alpha
u)}\biggr)
\]
for every $u\in(0,\uppi)$ shows that
\begin{eqnarray*}
\biggl(\frac{A_\beta}{A_\alpha}\biggr)' (\uppi z) & = &
\frac{1}{A_\alpha^2}\bigl((1-\beta^2)
A_\alpha- (1-\alpha^2)A_\beta+ A_\alpha A_\beta(A_\alpha- A_\beta)\bigr) (\uppi z)\\
& = & \frac{8 (\beta^2 - \alpha^2)(1-\alpha^2)(1-\beta^2) z^3}{\uppi^3
A_\alpha^2 (\uppi z)} \bigl(S_\alpha(z) S_\beta(z) S (z) - \uppi^2
S_{\alpha, \beta}(z)/4\bigr)
\end{eqnarray*}
for every $z\in(0,1),$ with the notations
\begin{eqnarray*}
S_\alpha(z) &=& \sum_{n\ge1}\frac{n^2}{(n^2 -\alpha^2z^2)(n^2 -
z^2)},\qquad S_\beta(z) = \sum_{n\ge1}\frac{n^2}{(n^2 -\beta^2z^2)(n^2
- z^2)},
\\
S(z) &=& \sum_{n\ge1}\frac{n^2}{(n^2 -\alpha^2z^2)(n^2 -
\beta^2z^2)},
\\
S_{\alpha, \beta}(z) &=& \sum_{n\ge1}\frac{n^2}{(n^2 -\alpha^2z^2)(n^2
- \beta^2z^2)(n^2 - z^2)} \cdot
\end{eqnarray*}
Since $S(z) \ge S(0) = \uppi^2/6,$ we see that (\ref{Final}) is true if
\[
S_\alpha(z) S_\beta(z) \ge3S_{\alpha, \beta} (z)/2
\]
for every $z\in(0,1).$ The latter follows for example, after isolating
the first term in each series and using the fact that $\uppi^2/6 \ge
3/2.$ We leave the details to the reader. Observe that the latter is
also true for $z = 0$ because $S_\alpha(0) = S_\beta(0) = \uppi^2/6$
and $S_{\alpha, \beta} (0) = \uppi^4/90.$
\end{pf}

\begin{pf*}{Proof of Theorem \ref{LCX1}}
We need to show that the density of
$\log(\ba^{-1}(U))$ is log-convex, in other words that the function
\[
\frac{x f_{\ba^{-1}(U)}'(x)}{f_{\ba^{-1}(U)}(x)}
\]
is increasing over its domain of definition which is
$[\alpha^\alpha\beta^\beta, +\infty)$. Since $(\log(\ba^{-1}))' = f$ is
a strictly absolutely monotonic function, the same holds for $\ba^{-1}$
and we set $\hba$ for its increasing reciprocal function. A computation
yields
\[
\frac{x f_{\ba^{-1}(U)}'(x)}{f_{\ba^{-1}(U)}(x)} = \frac{x \hba
''(x)}{\hba'(x)} = - 2 +
\biggl(\frac{\ba\ba''}{(\ba')^2}\biggr)(\hba(x)),
\]
which is an increasing function by Lemma \ref{trigo}.
\end{pf*}

\subsection{\texorpdfstring{Further properties of $\Wa$ and $\Va$}
{Further properties of W alpha and V alpha}}\label{s3.2}

In \cite{PSS}, the infinite divisibility of $\Wa$ is proved together
with a closed formula for its L\'evy measure. In this paragraph, we use
the latter expression to show that $\Wa$ is actually a translated ME.

\begin{THEO}
\label{CM1} The density of $\Wa$ is  CM.
\end{THEO}

As a consequence, we obtain the following reinforcement of Corollary
\ref{LCX2}, which entails that $\Za^\gamma$ is the product of a HCM
random variable and a positively translated ME random variable for
every $\gamma\ge\alpha/\beta$:

\begin{CORR}
\label{CM2} For every $s > 0,$ there exists $c_{\alpha, s} > 0$ such
that $\Va^s - c_{\alpha,s}$ is ME.
\end{CORR}

\begin{pf} Suppose first that $s = 1.$ Theorem \ref{CM1} shows that $\ga
(t) = f_{\Va}(\mathrm{e}^t)$ is CM, with the notation of Corollary
\ref{LCX2}. Since $\log x$ has a CM derivative, the classical Criterion
2 in \cite{F}, page 417, entails that $f_{\Va}(x) = \ga(\log x)$ is CM
over its domain of definition. It is clear from the definition of $\ba$
that the latter is $[\beta\alpha^{\alpha/\beta}, +\infty),$ so that
$\Va- \beta\alpha^{\alpha/\beta }$ is ME, by Proposition 51.8 in
\cite{S}. The general case $s > 0$ follows analogously in considering
the density of $s\Wa$ instead of $\Wa.$
\end{pf}

The proof of Theorem \ref{CM1} relies on the following lemma.

\begin{lemm}
\label{theetete} The L\'evy measure of $\Wa$ has a CM density given by
\[
w_\alpha(x) = \beta\int_1^\infty \mathrm{e}^{-\beta tx} ([t] - [\alpha
t] - [\beta t]) \,\mathrm{d}t,\qquad x\ge0,
\]
where $[t]$ stands for the integer part of any $t \ge0$.
\end{lemm}

\begin{pf}
From (\ref{GGCKant}) and (3.5), (3.6) and (3.7) in \cite{PSS} -- beware
our notation $\beta= 1- \alpha$, for every $\lbd\ge0$ one has
\begin{eqnarray*}
\EE[\mathrm{e}^{-\beta\lbd\Wa}] & = & \EE[\Za^{-\lbd\alpha}]/\EE[L^{
\lbd\beta
}]\\
& = & \frac{\Gamma(1+\lbd)}{\Gamma(1 +\lbd\alpha)\Gamma(1 +\lbd\beta
)}\\
& = & \exp-\biggl[a_\alpha\lbd+ \int_0^\infty(1- \mathrm{e}^{-\lbd x} -
\lbd x) \ga(x) \frac{\mathrm{d}x}{x}\biggr]
\end{eqnarray*}
for some constant $a_\alpha\in\rl,$ where
\[
\ga(x) = \frac{\mathrm{e}^{-x}}{1-\mathrm{e}^{-x}} -
\frac{\mathrm{e}^{-x/\alpha}}{1-\mathrm{e}^{-x/\alpha}} -
\frac{\mathrm{e}^{-x/\beta}}{1-\mathrm{e}^{-x/\beta}}
\]
is a nonnegative function -- see Lemma 3 in \cite{PSS} -- which is also
integrable with
\[
\int_0^\infty\ga(x) \,\mathrm{d}x = - (\alpha\log\alpha+
\beta\log\beta).
\]
Besides, $\ga(x)\to1$ as $x\to0$ so that one can rewrite
\begin{equation}
\label{LevyW} \EE[\mathrm{e}^{-\beta\lbd\Wa}] = \exp-\biggl[\hta\lbd+
\int_0^\infty(1- \mathrm{e}^{-\lbd x}) \ga(x)
\frac{\mathrm{d}x}{x}\biggr]
\end{equation}
for every $\lbd\ge0,$ where $\hta= \alpha\log\alpha+ \beta\log\beta$ is
the left-extremity of the support of the variable $\beta\Wa$ (this
shows that the above constant $a_\alpha$ is actually zero). Observe
that
\[
\ga(x) = \sum_{k=0}^\infty\bigl(\mathrm{e}^{-(k+1)x} -
\mathrm{e}^{-(k+1)x/\alpha} - \mathrm{e}^{-(k+1)x/\beta}\bigr) =
\int_1^\infty \mathrm{e}^{-tx} \mu_\alpha(\mathrm{d}t),
\]
where
\[
\mu_\alpha= \sum_{k=1}^\infty(\delta_k - \delta_{k/\alpha} - \delta
_{k/\beta})
\]
is a signed Radon measure over $\rl^+.$ Integrating by parts, we get
\[
\ga(x) = x\int_1^\infty \mathrm{e}^{-tx} \mu_\alpha([1,t))\,\mathrm{d}t
= x\int_1^\infty \mathrm{e}^{-tx} ([t] - [\alpha t] - [\beta t])
\,\mathrm{d}t.
\]
Putting everything together shows that the density of the L\'evy
measure of $\Wa$ is given by
\[
w_\alpha(x) = \beta\int_1^\infty \mathrm{e}^{-\beta tx} ([t] - [\alpha
t] - [\beta t])\,\mathrm{d}t,\qquad x\ge0.
\]
\upqed\end{pf}

\begin{pf*}{Proof of Theorem \ref{CM1}} From Lemma \ref{theetete} and Theorem
51.10 in \cite{S}, we already know that $\Wa- \hta$ belongs to the
class B, which is the closure of ME for weak convergence and
convolution. Moreover, one has
\[
w_\alpha(x) = \int_0^\infty \mathrm{e}^{-tx}
\theta_\alpha(t)\,\mathrm{d}t,\qquad x\ge0,
\]
with the notation $\theta_\alpha(t) = ([t/\beta] - [t] -
[(\alpha/\beta)t])\Un_{\{t\ge\beta\}},$ and it is clear that
\[
\int_0^1 \theta_\alpha(t)\frac{\mathrm{d}t}{t} <
\infty\quad\mbox{and}\quad 0 \le\theta _\alpha (t) \le1,\qquad t \ge0.
\]
From Theorem 51.12 in \cite{S}, this shows that $\Wa- \hta$ belongs to
the class ME itself, as required.
\end{pf*}

\begin{REMS}
\label{rkcm}  (a) In the terminology of Schilling, Song and
Vondra\v{c}ek \cite{SSV}, Lemma \ref{theetete} shows that the L\'evy
exponent of the ID random variable $\Wa$ is, up to translation, a
complete Bernstein function. Using Theorems 6.10 and 7.3 in \cite{SSV}
and simple transformations leads to the same conclusion that the
density of $\Wa$ is CM.

(b) Corollary \ref{CM2} and Theorem 51.10 in \cite{S} show that for
every $s > 0$ the L\'evy measure of the random variable $\Va^s$ has a
density $v_{\alpha,s}$ which is CM. We believe that $x\mapsto x
v_{\alpha,s}(x)$ is also CM, in other words, that $\Va^s$ is actually
GGC for every $s > 0.$ See Conjecture \ref{c1} below.

(c) The above Radon measure $\mu_\alpha$ is signed and not everywhere
positive, and we see from (\ref{LevyW}) that $\Wa$ is \textit{not} the
translation of a GGC random variable. This latter property might have
been helpful for a better understanding of the random variables
$\Va^s,$ although it is still a conjecture -- see Comment (1), page,
101 in \cite{B} -- that the transformation $x\mapsto \mathrm{e}^x - 1$
leaves the GGC property invariant.
\end{REMS}

\subsection{\texorpdfstring{The case $\alpha=1/2$ and three conjectures}
{The case alpha=1/2 and three conjectures}}\label{s3.3}

Taking $\alpha= 1/2$ in (\ref{GGCKant}) yields the factorizations
\[
Z_{1/2}^s = L^{-s} \times V_{1/2}^s
\]
for every $s\in\rl,$ where $V_{1/2}$ has the explicit density
\[
f_{V_{1/2}}(x) = \frac{1}{2\uppi x\sqrt{x-1/4}}\Un_{\{ x > 1/4\}}.
\]
Since $\log V_{1/2}$ has a log-convex density, the random variables
$V_{1/2}^s$ are not HM (this fact was already noticed in \cite{K}, see
the remark before Theorem 4.1 therein) and in particular not HCM.
However, the following proposition shows that $V_{1/2}^s$ is GGC at
least for $s\in[1/2, 1].$ We set $Y_s = V_{1/2}^s - 4^{-s}$ for every
$s > 0.$

\begin{prop}
\label{HCM12} The random variable $Y_s$ is HCM if and only if
$s\in[1/2, 1].$
\end{prop}

\begin{pf} Changing the variable and using the fact that every function $f$
on $\rl^+$ is HCM iff $x^\lbd f(1/x)$ is HCM for every $\lbd\in\rl,$
one sees that the following equivalences hold:
\[
Y_s \mbox{ is }\mathrm{HCM}\quad
\Leftrightarrow\quad\frac{1}{(x+1)\sqrt{(x+1)^t - 1}} \mbox{ is
}\mathrm{HCM} \quad\Leftrightarrow\quad\frac{1}{(x+1)\sqrt{(x+1)^t -
x^t}} \mbox{ is }\mathrm{HCM},
\]
with the notation $t = 1/s.$ Using the notation
\[
f_t (x) = \frac{1}{(x+1)\sqrt{(x+1)^t - x^t}}
\]
for every $t > 0,$ it is obvious that $f_1$ and $f_2$ are HCM. If now $
1< t< 2,$ then $x\mapsto(x+1)^t - x^t$ is obviously a Bernstein
function (i.e., a positive function with CM derivative), Criterion 2 in
\cite{F}, page 417, entails that $1/\sqrt{(x+1)^t - x^t}$ is CM and
$f_t$ is also CM. Since $f_t(0) = 1,$ this shows that $f_t$ is the
Laplace transform of some positive random variable and, by Theorem
5.4.1 in \cite{B}, $f_t$ will be HCM iff the latter variable is GGC. By
the Pick function characterization given in \cite{B}, Theorem 3.1.2,
setting $g_t(z) = f_t(-z)$ we need to show that $g$ is analytic and
zero-free on $\CC\setminus[0,\infty)$ and $\operatorname{Im}
(g_t'(z)/g_t(z))\ge0$ for $\operatorname{Im} z
> 0$ (in other words, that $g_t'/g_t$ is a Pick function -- see Section~2.4 in \cite{B}). The first point amounts to show that $ 1 -
(z/(z-1))^t$ does not vanish on $\CC\setminus[0,\infty),$ which is true
because $1 < t < 2.$ For the second point, we write
\[
\frac{g_t'(z)}{g_t(z)} = \frac{1}{1-z} + \frac{t}{2}
\biggl(\frac{(1-z)^u - (-z)^u}{(1-z)^{u+1} - (-z)^{u+1}}\biggr)
\]
with $u = t- 1 \in(0,1).$ Since $1/(1-z)$ is obviously Pick, we need
to show that
\[
h_u (z) = \bigl((1-z)^u - (-z)^u\bigr)\bigl((1-\bz)^{u+1} -
(-\bz)^{u+1}\bigr)
\]
is also Pick for every $u\in(0,1).$ We compute
\begin{eqnarray*}
\operatorname{Im} (h_u (z)) & = & \operatorname{Im} (z)
(\vert1-z\vert^{2u} + \vert z\vert^{2u}) - \operatorname{Im}
\bigl((-z)^u (1-\bz)^{u+1} + (-\bz)^{u+1}(1-z)^u\bigr)
\\
& = & \operatorname{Im} (z) \bigl(\vert(1-z)^{u} -
(-\bz)^{u}\vert^2\bigr) - \operatorname{Im}
\bigl((-z)^u (1-\bz)^u\bigr)\\
& = & \operatorname{Im} (z) \bigl(\vert(1-z)^{u} -
(-\bz)^{u}\vert^2\bigr) + \operatorname{Im} \bigl((\vert z\vert^2
-\bz)^u\bigr)
\end{eqnarray*}
which is nonnegative when $\operatorname{Im} z > 0$ because
$u\in(0,1).$ This shows the required property and proves that $Y_s$ is
HCM if $s\in[1/2, 1].$ Suppose now that $s > 1$ viz. $t < 1.$ If
$1/(x+1)\sqrt{(1+x)^t -1}$ were HCM, then it would also be HM and the
function
\[
y \mapsto\log(1 + \mathrm{e}^y) + \frac{1}{2} \log\bigl((1+
\mathrm{e}^y)^t - 1\bigr)
\]
would be convex on $\rl.$ Differentiating the above entails that the
function
\[
x \mapsto\frac{1}{x} + \frac{t}{2}\biggl(\frac{x^{t-1} - 1}{x^t
-1}\biggr) \sim-\frac{t}{2x^t}\qquad\mbox{as $x\to\infty$}
\]
would be nonincreasing on $[1,+\infty),$ a contradiction. Last, in view
of (ix), page 6 in \cite{B} it is easy to see that $Y_s$ is not HCM if
$s < 1/2,$ since $(1+z)^t - 1$ vanishes at least twice on
$\CC\setminus(-\infty,0]$ when $t > 2.$ This completes the proof.
\end{pf}

\begin{REMS}
(a) Except in the trivial cases $t = 1$ and $t=2,$ we
could not find any series of functions of the type described in (\ref
{dens}) converging pointwise to $f_t$ when $t\in[1,2].$

(b) From the above proposition, one might wonder if Kanter's variable
$\ba(U)^{-1/\alpha}$ is not a translated HCM in general. If it were
true, then (\ref{Kant}) would show that $\Za$ would be the product of a
HCM and a translated HCM for every $\alpha\le1/2$, so that from Theorem
6.2.2. in \cite{B} we would be quite close to the solution to
Bondesson's hypothesis. Nevertheless, to show the above property raises
computational difficulties significantly greater than those in Lemma \ref{trigo},
and it does not seem that this approach could be simpler than the one
suggested in \cite{B}, pages 88--89. See also Conjecture \ref{c1}
below.
\end{REMS}

The following corollary shows that there are GGC random variables
related to positive $\alpha$-stable powers with $\alpha\ge1/2.$

\begin{CORR}
\label{Anecdote} For every $\alpha\in[1/2, 1),$ the random variable
$(\Za\times Z_{1/2})^\alpha$ is GGC.
\end{CORR}

\begin{pf} Kanter's and Shanbhag--Sreehari's factorizations entail
\[
(\Za\times Z_{1/2})^\alpha\elaw(\Za/L)^\alpha\times V_{1/2}^\alpha\elaw
L^{-1} \times(Y_\alpha+ 4^{-\alpha})
\]
and from Proposition \ref{HCM12} and Theorem 4.3.1 in \cite{B}, we know
that $Y_\alpha+ 4^{-\alpha}$ is GGC because $\alpha\in[1/2, 1).$ Since
$L^{-1}$ is HCM, Theorem 6.2.1 in \cite{B} shows that $(\Za\times
Z_{1/2})^\alpha$ is GGC as well.
\end{pf}

As just mentioned, Proposition \ref{HCM12} shows that $V_{1/2}^s$ is
GGC for every $s\in[1/2, 1].$ From the conjecture made in \cite{B09},
Remark 3(ii), one may ask if this property does not remain true for
every $s > 1.$ Considering the variable $Y_s$ instead, this amounts to
show that the non-HCM function
\[
x \mapsto\frac{t}{\uppi(x+1)\sqrt{(x+1)^t - 1}}
\]
is still a GGC density for every $t\in(0,1).$ From Example 15.2.2. in
\cite{SSV}, one sees that the Laplace transform of the renewal measure
of a tempered $t$-stable subordinator subordinated through a
(1$/$2)-stable subordinator, which is given by
\[
x \mapsto\frac{1}{\sqrt{(x+1)^t - 1}},
\]
is a factor in this density. However, we could not find any convenient
expression for the Stieltjes (i.e., double Laplace) transform of the
above renewal measure which would entail that the Laplace transform of
$Y_s$ is HCM. If this renewal measure had a log-concave density, we
could apply Theorem~4.2.1 in \cite{B}, but this property does not seem
to be true even in the semi-explicit case $t = 1/2$ (inverse Gaussian
distribution).

From the above Remark \ref{rkcm}(b), we know that for every $s > 0$ and
every $\alpha\in(0,1)$ the L\'evy measure of $\Va^s$ has a density
$v_{\alpha,s}$ which is CM (more precisely, of the form given in
Theorem~51.12 of \cite{S}). Proposition \ref{HCM12} yields the
reinforcement that $x\mapsto x v_{1/2,s}(x)$ is CM for every $s\in[1/2,
1].$ Even though this is a very particular case, one might wonder if it
is not true in general.

\begin{CONJ}\label{c1} The random variable $\Va^s$ is GGC for every $s > 0$
and every $\alpha\in(0,1)$.
\end{CONJ}

When $\alpha= 1/2,$ a part of this conjecture can be rephrased in terms
of a more general question on Beta variables. Since
\[
(\operatorname{Beta} (\alpha_1, \alpha_2))^{-1} \elaw1 +
\frac{\operatorname{Gamma} (\alpha _2,1)}{\operatorname{Gamma}
(\alpha_1,1)}
\]
for every $\alpha_1, \alpha_2 > 0$ -- see, for exmple, \cite{B}, page
13, the HCM property for Gamma variables and Theorem 5.1.1 in \cite{B}
show that $(\operatorname{Beta} (\alpha_1, \alpha_2))^{-1}$ is always a
GGC. A particular case of the conjecture made in \cite{B09}, Remark
3(ii), is hence the following.

\begin{CONJ}\label{c2} The random variable $(\operatorname{Beta} (\alpha_1, \alpha_2))^{-s}$ is
GGC for every $\alpha_1, \alpha_2 > 0$ and every $s\ge1$.
\end{CONJ}

We could not find in the literature any result on the infinite
divisibility of negative powers of Beta variables. Proposition
\ref{HCM12} shows that $(\operatorname{Beta} (1/2, 1/2))^{-s}$ is a
positively translated HCM, hence GGC and ID, when $s\in[1/2, 1]$. The
same conclusion holds for other values of the parameters, not covering
the full range $\alpha_1, \alpha_2 > 0.$ From the identities
(\ref{GGCKant}) and Theorem 6.2.1 in \cite{B}, a positive answer to
Conjecture \ref{c1} would also show that $\Za^s$ is GGC for every
$s\ge\alpha/(1-\alpha).$ In view of our results in the previous section
for negative stable powers, it is tantalizing to raise the following general
conjecture.

\begin{CONJ}\label{c3} The random variable $\Za^\gamma$ is GGC for every
$\alpha \in(0,1)$ and \mbox{$\vert\gamma\vert\ge\alpha/(1-\alpha).$}
\end{CONJ}

If Bondesson's HCM conjecture is true, then $\Za^s$ is GGC for every
$\alpha\le1/2$ if $\vert s\vert\ge1.$ The above statement is stronger
since it takes values of $\alpha$ which are greater than 1$/$2 in
consideration, and since $\alpha/(1-\alpha) < 1$ when $\alpha< 1/2.$
Notice that if Conjecture \ref{c3} is true, then from Proposition
\ref{basic} we would also have $\Za^{-\gamma} \mbox{is GGC}
\Leftrightarrow\Za^{-\gamma} \mbox{is ID}$ for every $\gamma> 0.$ Last,
since the main theorem of \cite{TSS} shows that all positive stable
powers are unimodal, it would also be interesting to know if
$\Za^{\gamma}$ is still a GGC when $\gamma\in(0, \alpha/(1-\alpha)).$

\section{Related HCM densities}\label{s4}

In this section, we study two families of random variables which are
related to Bondesson's hypothesis, and have their independent interest.
For every $\alpha\in\rl,$ introduce the function
\[
\ga(u) = \frac{1}{u^{2\alpha} + 2 \cos(\uppi\alpha) u^\alpha+ 1}
\]
from $\rl^+$ to $\rl^+.$ The following elementary result might be well
known, although we could not trace any reference in the literature. As
for Proposition \ref{HCM12}, we could not find any constructive
argument either.

\begin{prop}
\label{Lamp} The function $\ga$ is HCM iff $\vert\alpha\vert\le 1/2.$
\end{prop}

\begin{pf}
Since $u^{2\alpha}\ga(u)$ HCM $\Leftrightarrow\ga(u)$ HCM and since
cosine is an even function, it is enough to consider the case
$\alpha\ge0.$ Suppose first $\alpha> 1/2.$ Rewriting
\[
\ga(u) = \frac{1}{(u^\alpha+
\mathrm{e}^{\mathrm{i}\uppi\alpha})(u^\alpha+
\mathrm{e}^{-\mathrm{i}\uppi \alpha})}
\]
shows that $\ga$ has two poles in $\CC/(-\infty,0]$ (because $\vert
(\alpha-1)/\alpha\vert< 1$) and from (ix), page 68 in \cite{B}, this
entails that $\ga$ is not HCM. Suppose next $0\le\alpha\le1/2.$ The
cases $\alpha= 0$ with $m_0(u) = 1/4$ and $\alpha= 1/2$ with
$m_{1/2}(u) = 1/(u+1)$ yield the HCM property explicitly, and we only
need to consider the case $0< \alpha< 1/2.$ Rewriting
\[
\ga(u) = \frac{u^{-\alpha}}{2 \cos(\uppi\alpha) + u^{-\alpha} +
u^\alpha}
\]
we see that it is enough to show that $\rho\dvtx  u\mapsto1/(c +
u^\alpha+ u^{-\alpha})$ is HCM for any $c > 0.$ Developing, one obtains
\[
\rho(uv)\rho(u/v) = \frac{1}{c^2 + u^{2\alpha} + u^{-2\alpha} +
v^{2\alpha} + v^{-2\alpha} + c(u^\alpha+ u^{-\alpha})(v^\alpha+
v^{-\alpha})}
\]
for any $u, v > 0.$ Since the function $v^{2\alpha} + v^{-2\alpha} +
c(u^\alpha+ u^{-\alpha})(v^\alpha+ v^{-\alpha})$ has CM derivative in
$w = v + 1/v$ for any fixed $u, c > 0$ (see \cite{BB}, page 183),
again Criterion 2 in \cite{F}, page 417, entails that the function
$\rho(uv)\rho(u/v)$ is CM in $w,$ so that $\rho$ is HCM.
\end{pf}

This proposition has several interesting consequences. Let us first
consider the random variable
\[
\Ya= \Ta^{1/\alpha}
\]
with the notation of Section~\ref{s2.2}, and recall that it is the
quotient of two independent copies of $\Za.$ The fact that $\Ya$ has a
closed density seems to have been first noticed in \cite{La}, in the
context of occupation time for certain stochastic processes. If
Bondesson's conjecture is true, then $\Ya$ is HCM whenever
$\alpha\le1/2,$ as a quotient of two independent HCM random variables
-- see Theorem 5.1.1 in \cite{B}. The next corollary shows that this is
indeed the case.

\begin{CORR}\label{co4.1}
The random variable $\Ya$ is HCM iff $\alpha\le1/2.$
\end{CORR}

\begin{pf} From a fractional moment identification, the density of $\Ya$ is
explicitly given by
\[
f_{\Ya}(x) = \frac{\sin\uppi\alpha x^{\alpha-1}}{\uppi(x^{2\alpha} + 2
x^\alpha\cos \uppi\alpha+ 1)} = \frac{\sin\uppi\alpha}{\uppi}
x^{\alpha-1}\ga(x)
\]
over $\rl^+$ -- see Exercise 4.21(3) in \cite{CY} already mentioned in
Section~\ref{s2}. The second derivative of $t \mapsto\log f_{\Ya}
(\mathrm{e}^t)$ equals
\[
\frac{-4\alpha^2(1 + \cos\uppi\alpha\cosh\alpha t)}{(\mathrm{e}^{\alpha
t} + 2 \cos \uppi\alpha+ \mathrm{e}^{-\alpha t})^2}
\]
and is not everywhere nonpositive whenever $\alpha
> 1/2,$ so that $\Ya$ is not HM, hence not HCM either. When $\alpha\le
1/2,$ the above Proposition \ref{Lamp} shows immediately that $f_{\Ya}$
is a HCM function, so that $\Ya$ is HCM as a random variable.
\end{pf}

We now turn our attention to the so-called Mittag--Leffler random
variables which were introduced in \cite{P}, and appeared since then in
a variety of contexts. Let
\[
\Ea(x) = \sum_{n=0}^{\infty} \frac{x^n}{\Gamma(1 +\alpha n)}
\]
be the classical Mittag--Leffler function with index $\alpha\in(0,1].$
The Mittag--Leffler random variable $\Ma$ has an explicit decreasing
density given by
\[
f_{\Ma}(x) = \alpha x^{\alpha-1} \Ea'(-x^\alpha)
\]
over $\rl^+.$ Its Laplace transform is also explicit: for every $\lbd
\ge0$ one has
\begin{equation}
\label{GGCMa} \EE[\mathrm{e}^{-\lbd\Ma}] = \frac{1}{1 + \lbd^\alpha} =
\exp \biggl[-\alpha\int_0^\infty(1 - \mathrm{e}^{-\lbd x})
\Ea(-x^\alpha) \frac {\mathrm{d}x}{x}\biggr]
\end{equation}
(see Remark 2.2 in \cite{P} and correct $x^k \to x^{\alpha k}$
therein). From the classical fact that $x \mapsto\Ea(-x)$ is the
Laplace transform of $\Za^{-\alpha}$ -- see, for example, Exercise
29.18 in \cite{S} -- and hence a CM function, this shows that $\Ma$ is
a GGC (with finite Thorin measure). This latter fact follows also from
the factorization
\[
\Ma\elaw\Za\times L^{1/\alpha},
\]
where $L \sim\operatorname{Exp} (1)$ (the latter is a direct
consequence of the first equality in (\ref{GGCMa}) -- see also the
final remark in \cite{P}), and from Theorem 6.2.1 in \cite{B} since
$\Za$ is GGC and $L^{1/\alpha}$ is HCM. Notice that in Example 3.2.4 of
\cite{BB}, the GGC property for $\Ma$ is also obtained in a slightly
more general context with the help of Pick functions.

\begin{CORR} The random variable $\Ma$ is HCM iff $\alpha\le1/2.$
\end{CORR}

\begin{pf} As a consequence of the above discussions, one has the classical
representation
\[
\Ea(-x^\alpha) = \EE[\mathrm{e}^{-x^\alpha\Za^{-\alpha}}]=
\EE[\mathrm{e}^{-x \Ya }] = \frac{\sin\uppi\alpha}{\uppi}\int_{\rl^+}
\frac{u^{\alpha-1} \mathrm{e}^{-xu} }{u^{2\alpha} +
2u^\alpha\cos\uppi\alpha+1} \,\mathrm{d}u
\]
so that the density of $\Ma$ writes
\[
f_{\Ma}(x) = \frac{\sin\uppi\alpha}{\uppi}\int_{\rl^+} \frac{u^\alpha
\mathrm{e}^{-xu} }{u^{2\alpha} + 2u^\alpha\cos\uppi\alpha+1}
\,\mathrm{d}u.
\]
From Proposition \ref{Lamp}, the function
\[
u \mapsto\frac{\sin\uppi\alpha u^\alpha}{\uppi(u^{2\alpha} +
2u^\alpha\cos\uppi\alpha +1)} = u f_{\Ya} (u)
\]
is HCM as soon as $\alpha\le1/2,$ so that it is also a widened GGC
density with the notations of Section~3.5 in \cite{BB}. By
Theorem~5.4.1 in \cite{B}, this shows that $f_{\Ma}$ is a HCM density
function whenever $\alpha\le1/2.$

There are two ways to prove that $\Ma$ is not HCM for $\alpha> 1/2.$
First, again from Theorem~5.4.1 in \cite{BB}, it suffices to show that
$u f_{\Ya} (u)$ is no more a widened GGC density when $\alpha> 1/2.$ If
it were true, then from Remark~6.1 in \cite{BB} the function $u^{1-
\delta} f_{\Ya} (u)$ would also be a widened GGC density for every
$\delta> 0.$ In particular, for every $\delta\in\,] \alpha, 1
+\alpha[$ the function $u^{1- \delta} f_{\Ya} (u)$ would be up to
normalization the density of a GGC. The derivative of the above
function equals
\[
\frac{- u^{\alpha- \delta- 1}((\delta+\alpha) u^{2\alpha} + 2\delta
u^\alpha\cos \uppi\alpha+ (\delta-\alpha))}{(u^{2\alpha} +
2u^\alpha\cos\uppi\alpha+1)^2}
\]
and is easily seen to vanish twice on $\rl^+$ if $\delta$ is close
enough to $\alpha> 1/2.$ This shows that the underlying variable is
bimodal and contradicts Theorem 52.1 in \cite{S} since all GGC's are
SD.

The second argument shows that $\Ma$ is not even HM when $\alpha> 1/2.$
From (7), page 207 in \cite{E} one has the asymptotic expansion
\begin{equation}
\label{Erde} \Ea(-z) = \sum_{n=1}^{N-1} \frac{(-1)^{n+1}
z^{-n}}{\Gamma(1-\alpha n)} + \mathrm{O} (z^{-N}), \qquad z \to+\infty,
\end{equation}
for any $N \ge2.$ Besides, by complete monotonicity, one can
differentiate this expansion term by term. Taking $N = 3,$ one obtains
\[
\bigl(z \Ea'''(-z) - \Ea''(-z)\bigr)\Ea'(-z) - z
\bigl(\Ea''(-z)\bigr)^2  =  \frac{-1}{\alpha^2 z^6
\Gamma(-\alpha)\Gamma(-2\alpha)} + \mathrm{O} (z^{-7})
\]
as $z\to\infty,$ and the leading term in the right-hand side is
positive when $\alpha> 1/2$: this shows that
$t\mapsto\Ea'(-\mathrm{e}^t)$ and, a fortiori, $t\mapsto\alpha
\mathrm{e}^{(\alpha-1)t}\Ea'(-\mathrm{e}^{\alpha t})$ are not
log-concave, so that $\Ma$ is not~HM.
\end{pf}

\begin{REMS}
(a) Since $\Za^{-\alpha}$ is not ID, the function $\Ea
(-x) = \EE[\mathrm{e}^{-x\Za^{-\alpha}}]$ is CM but never HCM. However,
repeating verbatim the above argument shows that
$x\mapsto\Ea(-x^\alpha)$ is HCM if and only if $\alpha\le1/2.$

(b) From the above proof, one has the equivalence
\[
\Ma\mbox{ \textit{is} }\mathrm{HM}
\quad\Longleftrightarrow\quad\Ma\mbox{ \textit{is} } \mathrm{HCM}
\quad\Longleftrightarrow\quad\alpha\le1/2.
\]
The variable $L^{1/\alpha}$ is always HCM, hence HM, and we know that
$\Za \mbox{ is HM} \Leftrightarrow\alpha\le1/2.$ Since the HM property
is closed by independent multiplication, this also proves that $\Ma$ is
HM as soon as $\alpha\le1/2.$ On the other hand, the influence of the
HM variable $L^{1/\alpha}$ in the product $L^{1/\alpha}\times\Za$ is
not important enough to make $\Ma$ HM when $\alpha> 1/2.$ From the
above equivalence, one can ask if the general identification
\[
\operatorname{HM} \cap\,\operatorname{GGC} = \operatorname{HCM}
\]
is true or not, and we could not find any counterexample. Such an
identification would show Bondesson's conjecture in view of the main
result of \cite{TS}.

(c) When $\alpha\le1/2$, the variable $\Ma^r$ is clearly HM, hence
unimodal, for every $r\neq0.$ On the other hand, when $\alpha> 1/2,$ it
is possible to find some $r\neq0$ such that $\Ma^r$ is not unimodal.
Indeed, it follows easily from the expansion (\ref{Erde}) that
\[
f_{\Ma^{-1/\alpha}}(0+) = \frac{-1}{\alpha\Gamma(-\alpha)},\qquad
f_{\Ma^{-1/\alpha }}'(0+) = \frac{-1}{\alpha\Gamma(-2\alpha)} > 0
\]
and $f_{\Ma^s}(0+) = +\infty$ $\forall s < -1/\alpha.$ By continuity of
$s \mapsto f_{\Ma^s}(x)$ for every $x > 0,$ this shows that $\Ma^s$ is
not unimodal for some $s < -1/\alpha$ close enough to $-1/\alpha.$
Hence, we have shown the further equivalence
\begin{equation}
\label{Equiva} \Ma\mbox{ \textit{is} }\mathrm{HM}\quad
\Longleftrightarrow\quad\Ma^r \mbox{\textit{ is unimodal for every}
}r\neq0.
\end{equation}
The same equivalence holds for $\Za$ as a consequence of the main
results of \cite{TS,TSS} and we can raise the following natural
question: \emph{For which class of positive random variables with
density does the equivalence \textup{(\ref{Equiva})} hold true?} It is
easy to see that if $X$ is such that $X^r$ is unimodal for every $r\neq
0,$ then $X$ is absolutely continuous.

(d) Reasoning exactly as at the beginning of Paragraph 2.2 in
\cite{TSS}, the decomposition
\[
\Ma^s = (L \times L^{\alpha-1})^{s/\alpha} \times b_\alpha^{-s/\alpha}
(U)
\]
and the concavity of $b_\alpha$ show that $\Ma^s$ is unimodal as soon
as $s \ge-\alpha,$ for every $\alpha\in(0,1).$ One might ask whether
$M_a^s$ is also unimodal for every $-1/\alpha\le s < -\alpha.$ Notice
that the above reasoning is not valid anymore, at least for $\alpha=
1/2$ because $b_{1/2}^r (U) = 2^r (\cos(U/2))^r$ is bimodal as soon as
$r > 1.$

(e) Both random variables $\Ya$ and $\Ma$ appear as special instances
of the Lamperti-type laws which were introduced in \cite{J}, to which
we refer for a thorough study. See also the numerous references
therein.
\end{REMS}

\section*{Acknowledgements}
This project was supported by King Saud University, Deanship of Scientific
Research,
College of Science Research Center.
Both authors wish to thank grants ANR-09-BLAN-0084-01 and CMCU-10G1503
for partial support.
The authors would like to thank L. Bondesson for the interest he
took in this paper, and also for providing them the manuscript~\cite{Bo99}.


%

\printhistory


\begin{thebibliography}{24}

\bibitem{B81}
%
\begin{barticle}[mr]
\bauthor{\bsnm{Bondesson},~\bfnm{Lennart}\binits{L.}} (\byear{1981}).
\btitle{Classes of infinitely divisible distributions and densities}.
\bjournal{Z.~Wahrsch. Verw. Gebiete} \bvolume{57} \bpages{39--71}.
\bnote{Corection and Addendum \textbf{59}, 277.}%
\bid{doi={10.1007/BF00533713}, issn={0044-3719}, mr={0623454}}%
\bptnote{check related}%
\bptok{imsref}%
\end{barticle}
%
\endbibitem

\bibitem{BB}
%
\begin{barticle}[mr]
\bauthor{\bsnm{Bondesson},~\bfnm{L.}\binits{L.}}
(\byear{1990}).
\btitle{Generalized gamma convolutions and complete monotonicity}.
\bjournal{Probab. Theory Related Fields}
\bvolume{85}
\bpages{181--194}.
\bid{doi={10.1007/BF01277981}, issn={0178-8051}, mr={1050743}}
\bptok{imsref}%
\end{barticle}
%
\endbibitem

\bibitem{B}
%
\begin{bbook}[mr]
\bauthor{\bsnm{Bondesson},~\bfnm{Lennart}\binits{L.}}
(\byear{1992}).
\btitle{Generalized Gamma Convolutions and Related Classes of
Distributions and
Densities}.
\bseries{Lecture Notes in Statistics}
\bvolume{76}.
\baddress{New York}: \bpublisher{Springer}.
\bid{doi={10.1007/978-1-4612-2948-3}, mr={1224674}}
\bptok{imsref}%
\end{bbook}
%
\endbibitem

\bibitem{Bo99}
%
\begin{bmisc}[auto]
\bauthor{\bsnm{Bondesson},~\bfnm{Lennart}\binits{L.}}
(\byear{1999}).
\bhowpublished{A problem concerning stable distributions. Technical report, Uppsala Univ.}
\bptok{imsref}%
\end{bmisc}
\endbibitem


\bibitem{B09}
%
\begin{barticle}[mr]
\bauthor{\bsnm{Bondesson},~\bfnm{Lennart}\binits{L.}}
(\byear{2009}).
\btitle{On univariate and bivariate generalized gamma convolutions}.
\bjournal{J.~Statist. Plann. Inference}
\bvolume{139}
\bpages{3759--3765}.
\bid{doi={10.1016/j.jspi.2009.05.015}, issn={0378-3758}, mr={2553760}}
\bptok{imsref}%
\end{barticle}
%
\endbibitem

\bibitem{CY}
%
\begin{bbook}[mr]
\bauthor{\bsnm{Chaumont},~\bfnm{L.}\binits{L.}} \AND
\bauthor{\bsnm{Yor},~\bfnm{M.}\binits{M.}}
(\byear{2003}).
\btitle{Exercises in Probability: A Guided Tour from Measure Theory to Random Processes, via
Conditioning}.
\bseries{Cambridge Series in Statistical and Probabilistic Mathematics}
\bvolume{13}.
\baddress{Cambridge}: \bpublisher{Cambridge Univ. Press}.
\bid{doi={10.1017/CBO9780511610813}, mr={2016344}}
\bptok{imsref}%
\end{bbook}
%
\endbibitem

\bibitem{CT}
%
\begin{barticle}[mr]
\bauthor{\bsnm{Cuculescu},~\bfnm{Ioan}\binits{I.}} \AND
\bauthor{\bsnm{Theodorescu},~\bfnm{Radu}\binits{R.}}
(\byear{1998}).
\btitle{Multiplicative strong unimodality}.
\bjournal{Aust. N. Z. J. Stat.}
\bvolume{40}
\bpages{205--214}.
\bid{doi={10.1111/1467-842X.00023}, issn={1369-1473}, mr={1631092}}
\bptok{imsref}%
\end{barticle}
%
\endbibitem

\bibitem{D}
%
\begin{barticle}[mr]
\bauthor{\bsnm{Demni},~\bfnm{Nizar}\binits{N.}}
(\byear{2011}).
\btitle{Kanter random variable and positive free stable distributions}.
\bjournal{Electron. Commun. Probab.}
\bvolume{16}
\bpages{137--149}.
\bid{doi={10.1214/ECP.v16-1608}, issn={1083-589X}, mr={2783335}}
\bptok{imsref}%
\end{barticle}
%
\endbibitem

\bibitem{De}
%
\begin{bmisc}[auto:STB|2012/06/08|12:49:54]
\bauthor{\bsnm{Devroye},~\bfnm{L.}\binits{L.}} (\byear{2009}).
\bhowpublished{Random variate generation for exponentially and
polynomially tilted stable distributions. \textit{ACM Transactions on Modeling and Computer Simulation}
\textbf{19} Article 18}.
\bptok{imsref}%
\end{bmisc}
%
\endbibitem

\bibitem{Di}
%
\begin{barticle}[mr]
\bauthor{\bsnm{Di{\'e}dhiou},~\bfnm{Alassane}\binits{A.}}
(\byear{1998}).
\btitle{On the self-decomposability of the half-{C}auchy distribution}.
\bjournal{J. Math. Anal. Appl.}
\bvolume{220}
\bpages{42--64}.
\bid{doi={10.1006/jmaa.1997.5790}, issn={0022-247X}, mr={1612063}}
\bptok{imsref}%
\end{barticle}
%
\endbibitem

\bibitem{E}
%
\begin{bbook}[auto:STB|2012/06/08|12:49:54]
\bauthor{\bsnm{Erd{\'e}lyi},~\bfnm{A.}\binits{A.}}
(\byear{1953}).
\btitle{Higher Transcendental Functions, Vol. III}.
\baddress{New York}: \bpublisher{McGraw-Hill}.
\bptok{imsref}%
\end{bbook}
%
\endbibitem

\bibitem{F}
%
\begin{bbook}[mr]
\bauthor{\bsnm{Feller},~\bfnm{William}\binits{W.}}
(\byear{1971}).
\btitle{An Introduction to Probability Theory and Its Applications, Vol.
II},
\bedition{2nd ed.}
\baddress{New York}: \bpublisher{Wiley}.
\bid{mr={0270403}}
\bptok{imsref}%
\end{bbook}
%
\endbibitem

\bibitem{J}
%
\begin{barticle}[mr]
\bauthor{\bsnm{James},~\bfnm{Lancelot~F.}\binits{L.F.}}
(\byear{2010}).
\btitle{Lamperti-type laws}.
\bjournal{Ann. Appl. Probab.}
\bvolume{20}
\bpages{1303--1340}.
\bid{doi={10.1214/09-AAP660}, issn={1050-5164}, mr={2676940}}
\bptok{imsref}%
\end{barticle}
%
\endbibitem

\bibitem{JRY}
%
\begin{barticle}[mr]
\bauthor{\bsnm{James},~\bfnm{Lancelot~F.}\binits{L.F.}},
\bauthor{\bsnm{Roynette},~\bfnm{Bernard}\binits{B.}} \AND
\bauthor{\bsnm{Yor},~\bfnm{Marc}\binits{M.}}
(\byear{2008}).
\btitle{Generalized gamma convolutions, {D}irichlet means, {T}horin measures,
with explicit examples}.
\bjournal{Probab. Surv.}
\bvolume{5}
\bpages{346--415}.
\bid{doi={10.1214/07-PS118}, issn={1549-5787}, mr={2476736}}
\bptok{imsref}%
\end{barticle}
%
\endbibitem

\bibitem{K}
%
\begin{barticle}[mr]
\bauthor{\bsnm{Kanter},~\bfnm{Marek}\binits{M.}}
(\byear{1975}).
\btitle{Stable densities under change of scale and total variation
inequalities}.
\bjournal{Ann. Probab.}
\bvolume{3}
\bpages{697--707}.
\bid{mr={0436265}}
\bptok{imsref}%
\end{barticle}
%
\endbibitem

\bibitem{La}
%
\begin{barticle}[mr]
\bauthor{\bsnm{Lamperti},~\bfnm{John}\binits{J.}}
(\byear{1958}).
\btitle{An occupation time theorem for a class of stochastic processes}.
\bjournal{Trans. Amer. Math. Soc.}
\bvolume{88}
\bpages{380--387}.
\bid{issn={0002-9947}, mr={0094863}}
\bptok{imsref}%
\end{barticle}
%
\endbibitem

\bibitem{P}
%
\begin{barticle}[mr]
\bauthor{\bsnm{Pillai},~\bfnm{R.~N.}\binits{R.N.}}
(\byear{1990}).
\btitle{On {M}ittag--{L}effler functions and related distributions}.
\bjournal{Ann. Inst. Statist. Math.}
\bvolume{42}
\bpages{157--161}.
\bid{doi={10.1007/BF00050786}, issn={0020-3157}, mr={1054728}}
\bptok{imsref}%
\end{barticle}
%
\endbibitem

\bibitem{RVY}
%
\begin{barticle}[mr]
\bauthor{\bsnm{Roynette},~\bfnm{B.}\binits{B.}},
\bauthor{\bsnm{Vallois},~\bfnm{P.}\binits{P.}} \AND
\bauthor{\bsnm{Yor},~\bfnm{M.}\binits{M.}}
(\byear{2009}).
\btitle{A family of generalized gamma convoluted variables}.
\bjournal{Probab. Math. Statist.}
\bvolume{29}
\bpages{181--204}.
\bid{issn={0208-4147}, mr={2792539}}
\bptok{imsref}%
\end{barticle}
%
\endbibitem

\bibitem{S}
%
\begin{bbook}[mr]
\bauthor{\bsnm{Sat\^{o}},~\bfnm{Ken-iti}\binits{K.}}
(\byear{1999}).
\btitle{L\'evy Processes and Infinitely Divisible Distributions}.
\bseries{Cambridge Studies in Advanced Mathematics}
\bvolume{68}.
\baddress{Cambridge}: \bpublisher{Cambridge Univ. Press}.
\bid{mr={1739520}}
\bptok{imsref}%
\end{bbook}
%
\endbibitem

\bibitem{SSV}
%
\begin{bbook}[mr]
\bauthor{\bsnm{Schilling},~\bfnm{Ren{\'e}~L.}\binits{R.L.}},
\bauthor{\bsnm{Song},~\bfnm{Renming}\binits{R.}} \AND
\bauthor{\bsnm{Vondra{\v{c}}ek},~\bfnm{Zoran}\binits{Z.}}
(\byear{2010}).
\btitle{Bernstein Functions: Theory and Applications}.
\bseries{De Gruyter Studies in Mathematics}
\bvolume{37}.
\baddress{Berlin}: \bpublisher{de Gruyter}.
\bid{mr={2598208}}
\bptnote{check year}
\bptok{imsref}%
\end{bbook}
%
\endbibitem

\bibitem{PSS}
%
\begin{barticle}[mr]
\bauthor{\bsnm{Shanbhag},~\bfnm{D.~N.}\binits{D.N.}},
\bauthor{\bsnm{Pestana},~\bfnm{D.}\binits{D.}} \AND
\bauthor{\bsnm{Sreehari},~\bfnm{M.}\binits{M.}}
(\byear{1977}).
\btitle{Some further results in infinite divisibility}.
\bjournal{Math. Proc. Cambridge Philos. Soc.}
\bvolume{82}
\bpages{289--295}.
\bid{issn={0305-0041}, mr={0448483}}
\bptok{imsref}%
\end{barticle}
%
\endbibitem

\bibitem{TS}
%
\begin{barticle}[mr]
\bauthor{\bsnm{Simon},~\bfnm{Thomas}\binits{T.}}
(\byear{2011}).
\btitle{Multiplicative strong unimodality for positive stable laws}.
\bjournal{Proc. Amer. Math. Soc.}
\bvolume{139}
\bpages{2587--2595}.
\bid{doi={10.1090/S0002-9939-2010-10697-4}, issn={0002-9939}, mr={2784828}}
\bptok{imsref}%
\end{barticle}
%
\endbibitem


\bibitem{TSS}
%
\begin{barticle}[auto:STB|2012/06/08|12:49:54]
\bauthor{\bsnm{Simon},~\bfnm{T.}\binits{T.}}
(\byear{2012}).
\btitle{On the unimodality of power transformations of positive
stable densities}.
\bjournal{Math. Nachr.}
\bvolume{285}
\bpages{497--506}.
\bptok{imsref}%
\bid{mr={2899640}}
\end{barticle}
%
\endbibitem


\bibitem{SVH}
%
\begin{bbook}[auto:STB|2012/06/08|12:49:54]
\bauthor{\bsnm{Steutel},~\bfnm{F.~W.}\binits{F.W.}} \AND
\bauthor{\bsnm{Van~Harn},~\bfnm{K.}\binits{K.}}
(\byear{2003}).
\btitle{Infinite Divisibility of Probability Distributions on the Real
Line}.
\baddress{New York}:
\bpublisher{Dekker}.
\bid{mr={2011862}}
\bptok{imsref}%
\end{bbook}
%
\endbibitem

\end{thebibliography}
\end{document}